\documentclass[a4paper, 12pt]{amsart}

\usepackage{amsmath} 
\usepackage{amssymb} 
\usepackage{amsthm}
\usepackage{graphicx}
\usepackage{tikz}
\usepackage{subcaption}    
   
    \usepackage{hyperref}
         
   \newtheorem{theorem}{Theorem}
   \newtheorem{lemma}[theorem]{Lemma}

\begin{document}
 \title[On partitions with $k$ corners]{On partitions with $k$ corners not containing the staircase with one more corner}

   \date{\today}
   \author{Emmanuel Briand}
   \address{Emmanuel Briand\\
     Departamento Matemática Aplicada I\\
     Escuela Técnica Superior de Ingeniería Informática\\
Avda. Reina Mercedes, S/N. 41012 Sevilla. SPAIN }
    \email{ebriand@us.es}
    \thanks{Funding: Partially supported by MTM2016-75024-P and FEDER, PID2020-117843GB-I00, and Proyectos I+D+i FEDER Andalucía US-1262169.}

\begin{abstract}
We give three proofs of the following result conjectured by Carriegos, De Castro-García and Muñoz Castañeda in their work on enumeration of control systems:   when \(\binom{k+1}{2} \le n < \binom{k+2}{2}\), there are as many partitions of $n$ with $k$ corners as  pairs of partitions \((\alpha, \beta)\) such
that \(\binom{k+1}{2} + |\alpha| + |\beta| = n\). 
\end{abstract}

    \maketitle

\section{Introduction}\label{introduction}

\begin{figure}[h]
  \centering
  \begin{subfigure}[t]{0.4\textwidth}
        \centering
        \includegraphics[scale=1.3]{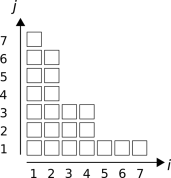}
        \caption{The diagram of the partition $\lambda = (7,4,4,2,2,2,1)$. The lengths of its rows are the parts $7$, $4$, $4$, $2$, $2$, $2$, $1$  of $\lambda$. This diagram has $4$ corners, whose positions $(i,j)$ are $(1, 7)$, $(2, 6)$, $(4, 3)$ and $(7, 1)$.}\label{fig1a}
    \end{subfigure}
    \qquad
    \begin{subfigure}[t]{0.4\textwidth}
           \centering
           \includegraphics[scale=1.3]{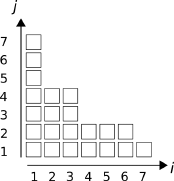}
           \caption{The diagram of the partition $\lambda'=(7,6,3,3,1,1,1)$. The parts of the $\lambda$ are the lengths of the columns in the diagram of $\lambda'$.}\label{fig1b}
    \end{subfigure}
    \caption{The diagrams of the partitions $\lambda=(7,4,4,2,2,2,1)$ and $\lambda'=(7,6,3,3,1,1,1)$, which are conjugate of each other. Each point $(i,j)$ of a diagram is represented as a square box centered at $(i,j)$.}\label{fig1}
  \end{figure}

\emph{Integer partitions} (finite weakly decreasing sequences of positive integers) are fundamental objects in enumerative combinatorics.
The terms of a partition are called its \emph{parts}.
Several parameters are attached to a partition, such as: its length, its weight (the sum of its parts), its largest part, the size of its Durfee square (i.e. the maximum $k$ such that the partition has $k$ parts $\ge k$) \ldots Another parameter is the number of distinct sizes of its parts. It is also the number of corners in the diagram of the partition. For instance,  the partition \((7,4,4,2,2,2,1)\) has parts of $4$ different sizes (7,4,2, and 1). Accordingly, its diagram has $4$ corners, as shown in Figure \ref{fig1a}. We will call for short \emph{partitions with $k$ corners} the partitions with parts of $k$ different sizes.

Recently, the problem of counting partitions with $k$ corners has arised in the context of the enumeration of linear control systems with coefficients in a commutative rings, in papers by Carriegos, De Castro-Garc\'{\i}a and Mu\~noz Casta\~neda \cite{Carriegos2016, Carriegos}. The present note is devoted to proving the following result that they conjectured.
\begin{theorem}[{\cite[Conjecture 30]{Carriegos}} ]\label{thm:conjectured}
\emph{When \(\binom{k+1}{2} \le n < \binom{k+2}{2}\), there are as many partitions of $n$ with $k$ corners as  pairs of partitions \((\alpha, \beta)\) such
that \(|\alpha| + |\beta| = n - \binom{k+1}{2}\).}
\end{theorem}

\begin{figure}[h]
  \centering
      \includegraphics[scale=0.8]{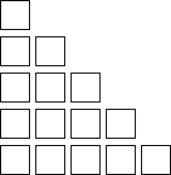}
        \caption{The diagram of $(5,4,3,2,1)$, which is the staircase partition of length $5$.}\label{fig:staircase}
  \end{figure}

Note that \(\binom{k+1}{2}\) is the size of the \emph{staircase
partition of length $k$}, which is \((k, k-1, k-2, \ldots, 1)\); see Figure \ref{fig:staircase}. 

We actually give three proofs of Theorem \ref{thm:conjectured}. The first one (Section \ref{proof-with-generating-functions}) is based on generating series. The second one (Section \ref{sec:Fine}) is based on an identity due to Fine on some statistics on partitions. The last one (Section \ref{bijective-proof}) is a bijective proof of the following more general result.

\begin{theorem}\label{thm:general}
For any $k \ge 0$ and $m \ge 0$, there are as many pairs of partitions $(\alpha, \beta)$ with $|\alpha|+|\beta|=m$ whose lengths fulfill $\ell(\alpha) + \ell(\beta) \le k$, as partitions of $m+\binom{k+1}{2}$ with $k$ corners whose diagrams do not contain the diagram of the staircase partition with length $k+1$. 
\end{theorem}

\begin{proof}[Proof of  Theorem \ref{thm:conjectured} from Theorem \ref{thm:general}] Consider $n < \binom{k+2}{2}$ and set $m = n - \binom{k+1}{2}$. Then $m \leq k$. Observe now that all pairs of partitions $(\alpha,\beta)$ with $|\alpha| + |\beta| = m$ fulfill $\ell(\alpha) + \ell(\beta) \le |\alpha| + |\beta| = m \le k$. Also, the diagram of a partition of $n$ cannot contain the diagram of the staircase partition of length $k+1$, since $n$ is smaller than $\binom{k+2}{2}=1+2+\cdots + (k+1)$. Theorem \ref{thm:conjectured} follows.\end{proof}

\section{Basic facts and notations}

\subsection{Partitions and their diagrams}\label{partitions-and-their-diagrams}

In this section, we recall classical operations and notations for integer partitions. See \cite{Pak}, \cite[I.1]{Macdonald}, \cite[p.58]{Stanley} or \cite[Ch. I]{AndrewsBook} for further details, basic definitions and general presentations.

Given a partition \(\lambda=(\lambda_1, \lambda_2, \ldots )\), its \emph{diagram}  is the set of integer
points \((i,j)\) such that \(1 \le j \le \ell(\lambda)\) and
\(1\le i \le \lambda_j\). In graphical representations, the points of the diagram of $\lambda$ are often drawn as square boxes centered at that point.
The \emph{conjugate} of $\lambda$ is the partition, denoted $\lambda'$, whose diagram is obtained from the diagram of $\lambda$ by applying the reflection that swaps the coordinates. The \emph{parts} of $\lambda$ are the nonzero terms $\lambda_1$, $\lambda_2$ \ldots and the \emph{length} of $\lambda$, denoted  with \(\ell(\lambda)\) , is the number of nonzero terms (non--necessarily distinct).
By $\lambda \vdash n$ we mean that $\lambda$ is a partition of $n$. We call $n$ the \emph{weight} of $\lambda$, and denote it with $|\lambda|$.

For instance, the partition $\lambda=(7,4,4,2,2,2,1)$ is a partition of $7+4+4+2+2+2+1=22$ (so $|\lambda|=22$ and $\lambda \vdash 22$), with length $7$ (denoted $\ell(\lambda)=7$) and parts $7$, $4$, $4$, $2$, $2$, $2$ and $1$. Its conjugate partition is $\lambda'=(7,6,3,3,1,1,1)$. The diagrams of $\lambda$ and $\lambda'$ are shown in Figure \ref{fig1}.

In sequences and partitions, we will make use of the notation $p^m$ for ``$m$ occurrences of $p$''. For instance, the partition $(7,4,4,2,2,2,1)$ will be also denoted  $(7, 4^2, 2^3, 1 )$.

Weakly decreasing sequences of non--negative integers with trailing zeroes will be identified with the partition obtained by deleting the trailing zeroes. For instance, $(7,4,4,2,2,2,1,0,0,0)$ will be identified to the partition $(7,4,4,2,2,2,1)$.

Consider two partitions $\alpha=(\alpha_1, \alpha_2, \ldots)$ and $\beta=(\beta_1, \beta_2, \ldots)$. Their \emph{sum} is the partition $\alpha + \beta = (\alpha_1+\beta_1, \alpha_2 + \beta_2, \ldots)$. For each $i$, let $a_i$ (resp. $b_i$) be the multiplicity of $i$ in $\alpha$ (resp. in $\beta$). The  \emph{union}
of the partitions $\alpha$ and $\beta$ is the partition (denoted $\alpha \cup \beta$) in which the multiplicity of $i$ is $a_i+b_i$. For instance, if $\alpha=(7,4,4,2,2,2,1)$ and $\beta=(4,2,1)$ then $\alpha+\beta=(7+4,4+2,4+1,2+0,2+0,2+0,1+0)=(11,6,5,2,2,2,1)$ and $\alpha \cup \beta=(7,4,4,4,2,2,2,2,1,1)$. The two operations are related by the identity $(\alpha+\beta)' = \alpha' \cup \beta'$.

\subsection{Generating series}\label{ref:basic_g_series}

Given any nonnegative integer $k$, the product
\[
\prod_{i=1}^{k} \frac{1}{1-q^i}
\]
expands as
\[
\sum_{a_1, a_2, \ldots, a_k} q^{1 a_1 + 2 a_2 + \cdots + k a_k},
\]
where the sum is carried over all $k$--tuples $(a_1, a_2, \ldots a_k)$ of nonnegative integers. Interpreting each integer $a_i$ as the multiplicity of $i$ as a part of a  partition with weight $1 a_1 + 2 a_2 + \cdots + k a_k$, we see that the series also writes $\sum_{\lambda} q^{|\lambda|}$, where the sum is over all partitions $\lambda$ whose parts are all $\le k$. It is thus the generating series for these partitions, according to their weight (i.e. for each $n$, the coefficient of $q^n$ in this series is the number of partitions of $n$ whose parts are all $\le k$).

Similarly, the infinite product
\[
\prod_{i=1}^{\infty} \frac{1}{1-q^i}
\]
is the generating series of all partitions, according to their weight.

\subsection{Corners of partitions}

A \emph{corner} of (the diagram of) \(\lambda\) is a
point \((i,j)\) in the diagram of \(\lambda\), such that neither
\((i+1,j)\) nor \((i, j+1)\) is in the diagram of \(\lambda\). The
partitions with parts of \(k\) distinct sizes  are exactly the partitions whose diagram has  \(k\) corners.

We will denote with \(\nu(n;k)\) the number of partitions of \(n\) with  \(k\) corners.

Let $\rho_k=(k, k-1, \ldots, 1)$. This partition is called the \emph{staircase partition of length  $k$}. Figure \ref{fig:staircase} shows the diagram of the staircase partition of length $5$. The staircase partition of length $k$ is the smallest partition with $k$ corners, in a sense made precise by the following lemma, that we will use implicitly in the sequel.

\begin{lemma}\label{lem:rhok}
  For any $k \ge 0$, the diagram of any partition with $k$ corners contains the diagram of $\rho_k$.
\end{lemma}

\begin{proof}
  Let $k \ge 0$ and let $\lambda$ be  a partition with $k$ corners. Then $\lambda $ is of the form $(p_1^{m_1} p_2^{m_2} \cdots p_k^{m_k}) $ with  $p_1 > p_2 > \cdots > p_k > 0$ and all $m_i > 0$.

  For each $i$, set $\mu_i = p_i - (k-i+1)$. Therefore, for all $i < k$, we have $\mu_i-\mu_{i+1}=p_i-p_{i+1}-1 \ge 0$. We have also $\mu_k=p_k-1 \ge 0$.
  This shows that $\mu = (\mu_1, \mu_2, \ldots, \mu_k)$ is a partition. We have:
  \[
  \lambda = (\rho_k + \mu) \cup (p_1^{m_1-1} p_2^{m_2-1} \cdots p_k^{m_k-1}).
  \]
This shows that the diagram of $\lambda$ contains the diagram of $\rho_k$.
  \end{proof}

We will also make repeated use of the following converse of Lemma \ref{lem:rhok}.
\begin{lemma}\label{lem:converse}
  Let $k \ge 0$, and let $\lambda$ be a partition. If the diagram of $\lambda$ does not contain the diagram of $\rho_{k+1}$, then $\lambda$ has at most $k$ corners. This is the case in particular for any partition $\lambda$ of weight less than $\binom{k+2}{2}$.
  \end{lemma}

\subsection{References on partitions with \texorpdfstring{$k$}{k} corners}

It has been remarked in \cite{Alladi} that the number of corners is a statistic on partitions that has been seldom considered. A notable and early exception is the work of MacMahon \cite{MacMahon} relating the enumeration of partitions according to their number of corners with number theory. MacMahon's work was expanded by Andrews \cite{Andrews} with a focus on asymptotics. MacMahon provided a formula for the generating series of partitions according to their weight and number of corners (identity \eqref{MacMahon gseries} below in Section \ref{proof-with-generating-functions}, where we make use of it). Refinements for this generating series were studied by Alladi \cite{Alladi}.  Let us mention also that part of MacMahon's results on this topic have been rediscovered by some authors unaware of his  work (for instance \cite{Kim}).

The sequence counting the partitions of $n$ with $k$ corners  is \href{https://oeis.org/search?q=A116608}{number A116608 in The On-Line Encyclopedia of Integer Sequences} \cite{oeis:triangle}.

Note finally that another possible name for the partitions with $k$ corners could have been ``partitions with $k$ distinct parts'', but this name is already taken and widely used for the partitions with $k$ distinct parts \emph{all of multiplicity $1$}, famous for having the same generating function as the partitions in odd parts.


\section{Proof with generating series}\label{proof-with-generating-functions}

In this section, we prove Theorem \ref{thm:conjectured} using the generating series of the numbers \(\nu(n; k)\), of partitions of $n$ with $k$ corners, which is defined as 
\[F(x,q) = \sum_{k, n} \nu(n; k)x^k q^n.\]
Our starting point will be the following expression for $F$, due to MacMahon.
\begin{theorem}[{\cite[p.90]{MacMahon}}; or see (3.2) in Andrews' account \cite{Andrews}]
  \label{MacMahon}
  The generating series $F$ is given by
  \begin{equation}\label{MacMahon gseries}
  F(x,q) =  \sum_{j=0}^{\infty} \frac{(x-1)^j q^{\binom{j+1}{2} }}{(1-q)(1-q^2)\cdots (1-q^j)} \cdot
\prod_{i=1}^{\infty} \frac{1}{1-q^i}\cdot
\end{equation}
\end{theorem}

\begin{proof}[Proof of Theorem \ref{thm:conjectured} from Theorem \ref{MacMahon}]
  Consider \(n\) and \(k\) as in Theorem \ref{thm:conjectured}. The number $\nu(n;k)$ is the coefficient of $x^n q^k$ in $F(x,q)$.

  In the right--hand side of \eqref{MacMahon gseries}, the summands with indices \(j < k\) have degree in \(x\) at most
\(k - 1\), and thus do not contribute to the coefficient of \(x^k q^n\) .
The expansions of the summands with indices \(j > k\) involve only monomials \(x^m q^d\)
with $d \ge \binom{k+2}{2} > n,$ and thus these summands do not contribute to the
coefficient of \(x^k q^n\).

Therefore, \(\nu(n;k)\) is the coefficient of \(x^k q^n\) in the summand
with index \(j = k\), which is 
\[\frac{(x-1)^k q^{\binom{k+1}{2}}}{(1-q)(1-q^2) \cdots (1-q^k)} \cdot \prod_{i=1}^{\infty} \frac{1}{1- q^i}\cdot\]
After expanding \((x-1)^k\) as $x^k + \text{terms of smaller degree in $x$}$ , we get that $\nu(n;k)$ is simply the coefficient
of \(q^n\) in \[
\frac{q^{\binom{k+1}{2}}}{(1-q)(1-q^2) \cdots (1-q^k)} \cdot \prod_{i=1}^{\infty} \frac{1}{1- q^i}
\] which is the coefficient of \(q^{n-\binom{k+1}{2}}\) in
\begin{equation}\label{series product}
\frac{1}{(1-q)(1-q^2) \cdots (1-q^k)} \cdot \prod_{i=1}^{\infty} \frac{1}{1- q^i}\cdot
\end{equation}
The left factor in \eqref{series product} is the generating series for the partitions $\lambda$ with parts $\le k$; the right factor is the generating series for all partitions $\mu$ (see section \ref{ref:basic_g_series}).
Therefore, for any $h$, the coefficient of \(q^h\) in \eqref{series product} is the number of pairs
of partitions \((\lambda, \mu)\) such that \(|\lambda|+|\mu| = h\), and
\(\lambda\) has all parts \(\le k\). For \(h \le k\), the condition on
the sizes of the parts of $\lambda$ can be dropped. This is the case in particular for
\(h = n - \binom{k+1}{2}\), since
\(n < \binom{k+2}{2} = \binom{k+1}{2} + k +1\).
\end{proof}


\section{Proof from statistics on partitions}\label{sec:Fine}

We now give another proof of Theorem \ref{thm:conjectured}, based on the following result due to  Fine. 
\begin{theorem}[{\cite[Theorem 4 in Chapter 2]{Fine}}]\label{thm:Fine}
For any $r \ge 0$,
\begin{equation}\label{eq:Fine}
\sum_{\lambda \vdash n} \binom{Q(\lambda)}{r} 
=  
\sum_{\lambda \vdash n} m_1(\lambda) m_2(\lambda) \cdots m_r(\lambda),
\end{equation}
where  $Q(\lambda)$ is the number of corners of $\lambda$, and $m_i(\lambda)$ stands for the multiplicity of $i$ as a part of $\lambda$.
\end{theorem}

\begin{proof}[Proof of Theorem \ref{thm:conjectured} from Theorem \ref{thm:Fine}] Consider $n$ and $k$ such that 
\[
\binom{k+1}{2} \le n < \binom{k+2}{2}.
\]
After Lemma \ref{lem:converse}, the diagram of any partition of $n$ has at most $k$ corners. Apply Theorem \ref{thm:Fine} with $r=k$. For any $\lambda \vdash n$, either $\lambda$ has less than $k$ corners, and then $\binom{Q(\lambda)}{k}=0$, or $\lambda$ has exactly $k$ corners, and then $\binom{Q(\lambda)}{k}=1$. The left-hand side in \eqref{eq:Fine} is thus $\nu(n; k)$. 

The right--hand side is
\begin{equation}\label{RHS}
\sum_{\lambda \vdash n} m_1(\lambda) m_2(\lambda) \cdots m_k(\lambda).
\end{equation}
Note that if a partition $\lambda \vdash n$ has some part $j>k$, then at least one of its parts $i \le k$ must have multiplicity $0$ in $\lambda$ (otherwise $n \le \binom{k+2}{2}$ cannot be fulfilled). Such a partition does not contribute to the sum \eqref{RHS}.

Consider the map $(\alpha, \beta) \mapsto \alpha \cup \beta \cup \rho_k$ from the pairs of partitions $(\alpha,\beta)$ with parts $\le k$. If $\alpha=(k^{a_k}\cdots 2^{a_2} 1^{a_1})$ and $\beta=(k^{b_k} \cdots 2^{b_2}1^{b_1})$, the image of $(\alpha,\beta)$ is the partition $(k^{a_k+b_k+1}\cdots 2^{a_2+b_2+1} 1^{a_1+b_1+1})$. This shows that the number of preimages of a partition $\lambda$, under this map,  is $m_1(\lambda) m_2(\lambda) \cdots m_k(\lambda)$. Indeed, each pair mapped to $\lambda$ is determined by the choice of the multiplicities $a_i$, that must fulfill $0 \le a_i \le m_i(\lambda)-1$ for all $i$, whence there are $m_i(\lambda)$ choices for $a_i$. The sum \eqref{RHS} now interprets as the cardinality of the inverse image of the set of all partitions $\lambda \le n$. This inverse image is the set of all pairs of partitions $(\alpha,\beta)$ with parts $\le k$ such that $|\alpha|+|\beta| +|\rho_k| = n$. This last condition simplifies as $|\alpha|+ |\beta|=n-\binom{k+1}{2}$. Finally, since $n < \binom{k+2}{2}$, we have that $n-\binom{k+1}{2} < k+1$. Therefore, the condition that the  parts of $\alpha$ and $\beta$ are all at most $k$ can be dropped, as it is already implied by the condition on their weights. The statement of Theorem \ref{thm:conjectured} is obtained.
\end{proof}


\section{Bijective proof}\label{bijective-proof}

It has been shown in the introduction that Theorem  \ref{thm:conjectured} follows from the more general Theorem \ref{thm:general}.
In this section, a bijective proof of Theorem \ref{thm:general} is provided.


We will use the \emph{border coordinates} for partitions, that we introduce now.

    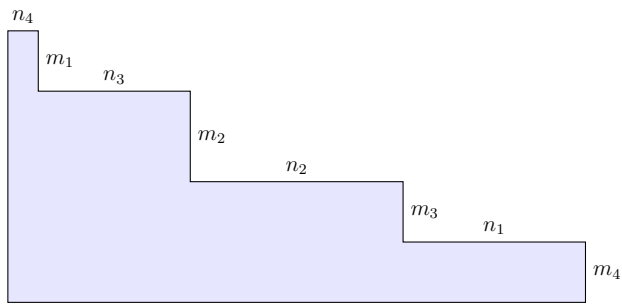
\begin{figure}
      \begin{center}
      \begin{tikzpicture}[scale=0.8]
\draw[fill=blue!10] (0,0) node (O) {}
-- ++(9.5,0) 
-- ++(0,1) node [midway, right, scale=0.7] () {$m_4$}
-- ++(-3,0) node [midway, above, scale=0.7] () {$n_1$} 
-- ++(0,1) node [midway, right, scale=0.7] () {$m_3$}
-- ++(-3.5,0) node [midway, above, scale=0.7] () {$n_2$} 
-- ++(0,1.5) node [midway, right, scale=0.7] () {$m_2$}
-- ++ (-2.5,0) node [midway, above, scale=0.7] () {$n_3$}
-- ++(0,1) node [midway, right, scale=0.7] () {$m_1$}
-- ++(-0.5,0) node [midway, above, scale=0.7] () {$n_4$}
-- cycle;
      \end{tikzpicture}
      \end{center}
      \caption{The diagram of a partition with $4$ corners, with its border coordinates $(m_1, m_2, m_3, m_4; n_1, n_2, n_3, n_4) $.}\label{fig border}
      \end{figure}
    
Let \(\lambda\) be a partition with \(k\) corners. Let \(p_1\),
\(p_2,\ldots, p_k\) be the distinct parts of \(\lambda\), ordered increasingly;
\(p_1 < p_2 < \cdots < p_k\). Let \(q_1 < q_2 < \cdots < q_k\) be the
distinct parts of the conjugate partition \(\lambda'\). For each \(i\),
let \(m_i\) (resp. \(n_i\)) be the multiplicity of \(p_i\) (resp.
\(q_i\)) in \(\lambda\) (resp. \(\lambda'\)). We call the pair of
sequences \((m_1, m_2, \ldots, m_k ; n_1, n_2,\ldots, n_k)\) the
\emph{border coordinates} of \(\lambda\), since they are the lengths of
the vertical and horizontal segments in the border of the diagram of
\(\lambda\) (see Figure \ref{fig border}).

The border coordinates $n_i$ are directly obtained as differences of consecutive parts. Precisely, 
\begin{equation}\label{border from parts}
n_k=p_1, \quad \text{ and for all } i > 1, \; n_{k+1-i} = p_{i}-p_{i-1}.
\end{equation}
Indeed, $p_1$, $p_2$, \ldots, $p_k$ are the first coordinates of the corners listed from left to right, while $n_{k-1}$, $n_{k-2}$, \ldots, $n_1$ are the differences between the first coordinates of consecutive corners. The same relation holds between the $m_i$ and the $q_i$.

\begin{lemma}\label{lem:border}
  Let \(\lambda\) be a partition with \(k\) corners and
border coordinates 
\[(m_1, m_2, \ldots, m_k; n_1, n_2, \ldots, n_k).\]
Let $p_1 < p_2 < \cdots < p_k$ be the parts of $\lambda$, and let $q_1 < q_2 < \cdots < q_k$ be the parts of $\lambda'$.
\begin{enumerate}
\item 
 Let \(\gamma\) be a partition whose parts
are all among the parts of \(\lambda\). For each \(i\), let \(c_i\) be the
multiplicity of \(p_i\) in \(\gamma\). Then the border coordinates of
\(\lambda \cup \gamma\) are
\[(m_1 + c_1, m_2 + c_2, \ldots, m_k + c_k; n_1, n_2, \ldots, n_k).\] 
\item 
Let \(\gamma\) be a partition such that all parts of \(\gamma'\) are
among the parts of \(\lambda'\). For each \(i\), let \(d_i\) be the
multiplicity of \(q_i\) in \(\gamma'\). Then the border coordinates of
\(\lambda + \gamma\) are
\[(m_1 , m_2 , \ldots, m_k ; n_1 + d_1, n_2 + d_2, \ldots, n_k + d_k).\]
\end{enumerate}
\end{lemma}

\begin{proof}
  We prove part 1. We have $\lambda = (p_k^{m_k} \cdots p_2^{m_2} p_1^{m_1})$ and $\gamma = (p_k^{c_k}\cdots  p_2^{c_2} p_1^{c_1})$. Therefore
  \[(\lambda \cup \gamma) = (p_k^{m_k+c_k} \cdots p_2^{m_2+c_2} p_1^{m_k+c_1}).
  \]
  This gives the first half of the border coordinates of $\lambda \cup \gamma$. For the second half, observe that, after \eqref{border from parts}, the $n_i$  do not depend on the multiplicities of the $p_i$, and thus remain unaffected by the union with $\gamma$.
  
  Part 2 is straightforwardly deduced from part 1 by means of  the identity 
  \(\lambda + \gamma = (\lambda' \cup \gamma')'\), and noting that conjugation swaps the  $m_i$'s and the $n_i$'s in the border coordinates.
\end{proof}

\begin{figure}
  \begin{center}
    \includegraphics[width=\textwidth]{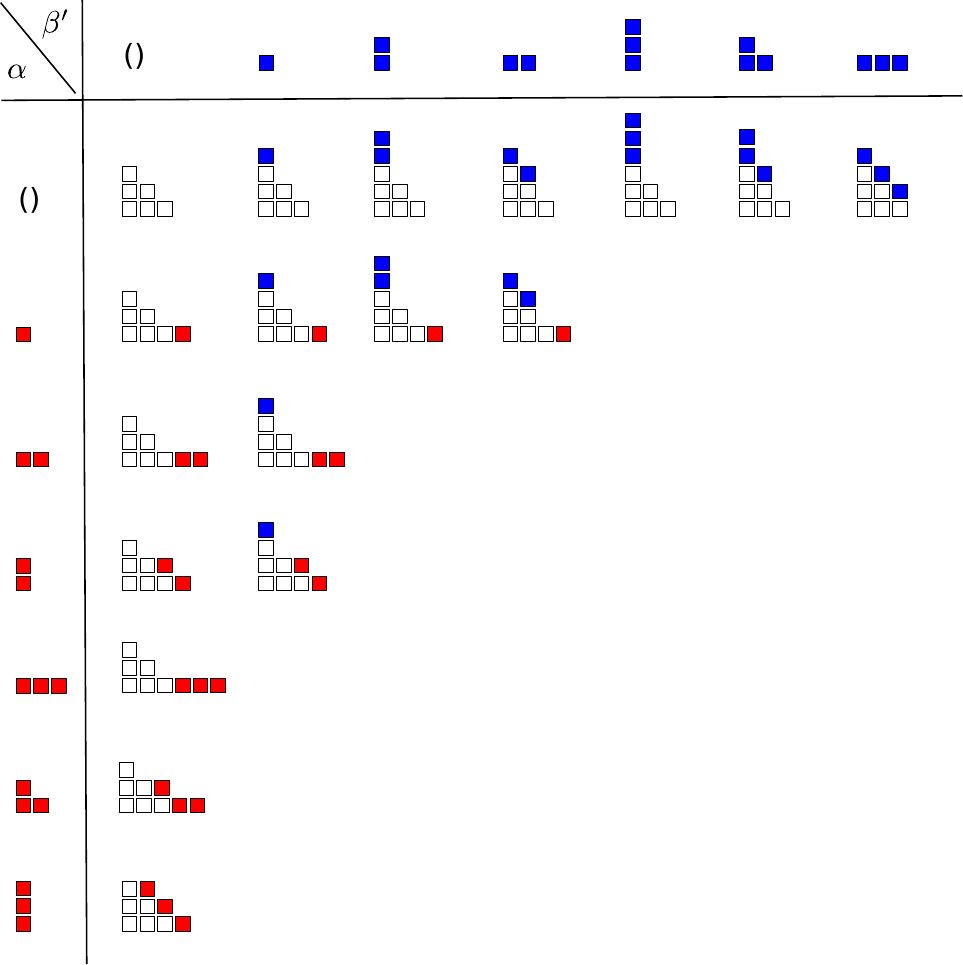}
    \end{center}
  \caption{The bijection $(\alpha,\beta) \mapsto (\rho_3 \cup \beta') +\alpha$ described in Lemma \ref{lem:general}. Its puts in correspondence the  pairs of partitions $(\alpha,\beta)$ with $|\alpha|+|\beta| \le 3$,  with the partitions $\lambda$ with $3$ corners and weight $< \binom{5}{2}$,  whose diagrams do not contain the diagram of $\rho_4$. This illustrates Theorem \ref{thm:conjectured} for $k=3$.
Rows correspond to the partitions $\alpha$, and columns to the partitions $\beta$. The conjugate partition $\beta'$ is drawn. The entry in row $\alpha$ and column $\beta$ is the diagram of the partition $(\rho_3 \cup \beta') + \alpha$. The staircase partition $\rho_3$ is left in white, and the new points of the diagram coming from $\alpha$ and $\beta'$ are in color.
  }\label{tabla}
 \end{figure}

Theorem \ref{thm:general} follows straightforwardly from the more precise lemma below.    
\begin{lemma}\label{lem:general}
  The map
\((\alpha, \beta) \mapsto (\rho_k \cup \beta') + \alpha\) establishes a
bijection between the pairs of partitions \((\alpha, \beta)\) such that
\(\ell(\alpha) + \ell(\beta) \le k\), and the partitions \(\lambda\)
with $k$ corners whose diagrams do not contain the diagram of
\(\rho_{k+1}\).
\end{lemma}
The case $k=3$ of the  bijection defined in Lemma \ref{lem:general} is shown in Figure \ref{tabla}. 

\begin{proof}    
Let \(\alpha\) and \(\beta\) be two partitions whose lengths have sum at
most \(k\). There exist \(p \ge\ell(\alpha)\) and \(q \ge \ell(\beta)\)
with \(p+q = k\). Set \(b_i\) (resp. \(a_i\)) for the multiplicity of
\(i\) in \(\beta'\) (resp. \(\alpha'\)).

Since $\beta$ has length at most $k$, all parts of $\beta'$ are smaller than or equal to $k$, and thus all parts of $\beta'$ are among those of $\rho_k$. Lemma \ref{lem:border} applies: the border coordinates of $\rho_k \cup \beta'$ are
$(1+b_1, 1+b_2, \ldots, 1 + b_k; 1^k)$.

The corners of the diagram of $\rho_k$ are the pairs $(i, k+1-i)$ for $i$ from $1$ to $k$.
Since $\beta$ has length at most $q$, i.e. $\beta'$ has all its parts smaller than or equal to $q$, performing the union with $\beta'$ does not affect the columns of the diagram whose indices $i$ fulfill $i > q$.
Therefore, all pairs $(i, k+1-i)$ for $i > q$ are still corners of $\rho_k \cup \beta'$.
As a consequence, the numbers $k+1-i$ for $i > q$, which are the numbers $j \le p$, are still parts of  $(\rho_k \cup \beta')'$.
This shows that all parts of $\alpha'$ are among the parts of $(\rho_k \cup \beta')'$, allowing to apply again Lemma \ref{lem:border}. The conclusion now is that the border coordinates of $(\rho_k \cup \beta') + \alpha$ are
\[
(1+b_1, 1+b_2, \ldots, 1 + b_k; 1+a_1, 1+a_2, \ldots, 1+a_k).
\]
This shows clearly that the map $\Phi$ that sends each pair of partitions $(\alpha,\beta)$ whose lengths have sum at most \(k\)  to $ (\rho_k + \beta')+ \alpha$ is injective.
Indeed, the multiplicities $a_i$ and $b_i$ of $\alpha'$ and $\beta'$  can be read from the border coordinates of $\Phi(\alpha,\beta)$.

Besides, since the union with $\beta'$ affects only the first $q$ columns of the diagram, and the sum with $\alpha$ affects only the first $p$ rows, the point $(q+1, p+1)$, that is not in the diagram of $\rho_k$, is still not in the diagram of $(\rho_k \cup \beta') + \alpha$.
Since $q+1+p+1=k+2$, this point $(q+1,p+1)$ lies in the diagram of $\rho_{k+1}$. This shows that the diagram of  $\Phi(\alpha,\beta)$ does not contain the diagram of $\rho_{k+1}$. On the other hand, $\Phi(\alpha,\beta)$ has $k$ corners exactly (since it has $2k$ border coordinates). We conclude that $\Phi$ takes its values in the set $\mathcal{S}_k$  of all  partitions with $k$ corners whose diagrams do not contain the diagram of $\rho_{k+1}$.

Let us now show that $\Phi$ has image $\mathcal{S}_k$ exactly. Let \(\lambda\) be a partition with \(k\) corners, whose
diagram does not contain the diagram of \(\rho_{k+1}\). 
There exists \((i_0, j_0)\), lying in 
the diagram of $\rho_{k+1}$, but not in the diagram of \(\lambda\). 
Since $\lambda$ has $k$ corners, its diagram contains the diagram of \(\rho_k\). Therefore  \((i_0, j_0)\) is not in the diagram of $\rho_k$, and belongs  to the set difference of the diagram of $\rho_{k+1}$ and the diagram of $\rho_k$. Thus \(i_0 + j_0 = k + 2\).
Moreover, any point $(i,j)$ of the diagram, and in particular any corner, must fulfill $i < i_0$ or $j < j_0$. Let \(p\) be
the number of corners \((i,j)\)  with
\(i < i_0\), and let \(q\) be the number of corners \((i,j)\) with
\(j < j_0\). Then $k \le p + q$. We have \(p \le i_0-1\) since there is at most one corner in each column. Similarly, \(q \le j_0-1\), because there is at most one
corner in each row.  Altogether, we get
\(k \le p + q \le (i_0-1) + (j_0-1) = k\). As a consequence,
\(k = p + q\) and \(p = i_0-1\), \(q = j_0 -1\). There is one corner
in each of the first \(p\) columns and one corner in each of the first $q$ rows, and no corner belongs at the same time to some of the first $p$ columns and to some of the first $q$ rows. We conclude that \(\lambda\) has
border coordinates of the form
\((m_1,\ldots,m_q, 1^p; n_1, \ldots, n_p, 1^q)\) for some positive
numbers \(m_i\) and \(n_i\). This coincides with the border coordinates of 
\((\rho_k \cup \beta')+\alpha\) for
\(\beta'=(q^{m_q-1} \cdots  2^{m_2-1} 1^{m_1-1}  )\) and
\(\alpha'=(p^{n_p-1}  \cdots 2^{n_2-1} 1^{n_1-1})\).
Clearly, partitions are determined by their border coordinates. Therefore, 
$\lambda$ is equal to $(\rho_k \cup \beta')+\alpha$ for the partitions $\alpha'$ and $\beta'$ defined above (and $\alpha=(\alpha')'$). This shows that all $\lambda\in \mathcal{S}_k$ are in the image of $\Phi$, which terminates the proof.
\end{proof}

\section*{Acknowledgments}

Thanks to the organizers of the \emph{XI Encuentro Andaluz de Matemática Discreta} held in Sevilla in February 2020, that allowed the author to hear from Professor Carriegos the conjecture considered in the present note. Thanks to Professor Carriegos for his interest in this work. Thanks to the anonymous reviewers for their valuable  comments.

\bibliographystyle{elsarticle-num}
\bibliography{corners}

\end{document}